\documentclass{article}%
\usepackage{amsfonts}
\usepackage{amsmath}
\usepackage{hyperref}
\usepackage{amssymb}
\usepackage{graphicx}%
\providecommand{\U}[1]{\protect\rule{.1in}{.1in}}
\begin{document}

\title{\vspace{-0.5in}Matrix exponentials, SU(N) group elements, and real polynomial roots}
\author{T. S. Van Kortryk\\120 Payne Street, Paris, MO65275\\{\small vankortryk@gmail.com}}
\date{}
\maketitle

\begin{abstract}
The exponential of an $N\times N$\ matrix can always be expressed as a matrix
polynomial of order $N-1$. \ In particular, a general group element for the
fundamental representation of $SU(N)$\ can be expressed as a matrix polynomial
of order $N-1$\ in a traceless $N\times N$\ hermitian generating matrix, with
polynomial coefficients consisting of elementary trigonometric functions
dependent on $N-2$ invariants in addition to the group parameter. \ These
invariants are just angles determined by the direction of a real $N$-vector
whose components are the eigenvalues of the hermitian matrix. \ Equivalently,
the eigenvalues are given by projecting the vertices of an $\left(
N-1\right)  $-simplex onto a particular axis passing through the center of the
simplex. \ The orientation of the simplex relative to this axis determines the
angular invariants and hence the real eigenvalues of the matrix.

\end{abstract}

\begin{quote}
{\small \textquotedblleft Let us revisit Euclid. \ Let us discover for
ourselves a few of the newer results. \ Perhaps we may be able to recapture
some of the wonder and awe that our first contact with geometry
aroused.\textquotedblright\ --- H S M Coxeter}
\end{quote}

Curtright and Zachos (CZ) wrote a brief summary \cite{CZ}\ of essential
elementary features for the triplet representation \cite{MGM} of $SU(3)$,
thereby distilling some older results \cite{Ilamed,MSW,Rosen}. \ Here, I show
how their main results may be extended to almost any $N\times N$ matrix,
thereby embellishing \cite{Ilamed,Kusnezov,Laufer}. \ A matrix polynomial form
for an exponentiated matrix can always be expressed succinctly in terms of a
\textquotedblleft response function\textquotedblright\ that encodes the
eigenvalues of the matrix.

In particular, I show how the CZ results may be extended to the fundamental
representation of $SU(N)$ for any $N$. \ I show that a polynomial form for any
group element in the fundamental representation can be expressed in terms of a
response function\ that encodes the real eigenvalues of the hermitian matrix
that generates the group element. \ In addition, I provide a clear geometrical
picture of the relevant group invariants in terms of the elementary properties
of an $\left(  N-1\right)  $-simplex \cite{Coxeter}, as a simple
generalization of Vi\`{e}te's venerable results for the real roots of a cubic
equation \cite{Viete}. \ I give specific results for $SU\left(  N\leq5\right)
$.

The exponential of an $N\times N$ matrix $M$ can be written as a matrix
polynomial as a consequence of the Cayley-Hamilton theorem \cite{Cayley}, or
more directly, as the result of the Lagrange-Sylvester projection matrix
method \cite{Sylvester}. \ A reasonably compact polynomial form is given by%
\begin{gather}
\exp\left(  itM\right)  =\sum_{n=0}^{N-1}M^{n}E_{n}\left(  t\right)
\ ,\label{ExponentialPoly}\\
E_{n}\left(  t\right)  =\sum_{m=0}^{N-1-n}\left(  -1\right)  ^{m}S_{m}\left[
\lambda\right]  \left(  -i\frac{d}{dt}\right)  ^{N-1-n-m}F\left(  t\right)
\ .
\end{gather}
The invariant functions appearing in the $E_{n}$ coefficients are given by
symmetric polynomials,%
\begin{equation}
S_{m}\left[  \lambda\right]  \equiv\sum_{1\leq k_{1}<k_{2}<\cdots<k_{m}\leq
N}\lambda_{k_{1}}\lambda_{k_{2}}\cdots\lambda_{k_{m}}\ ,
\end{equation}
where $\lambda_{k}$ for $k=1,\cdots,N$ are the eigenvalues of $M$, and by
various derivatives of the response function for the matrix, $F\left(
t\right)  $. \ The latter is defined in terms of the characteristic function
$C\left(  z\right)  $ as \
\begin{equation}
F\left(  t\right)  =\sum_{k=1}^{N}\frac{\exp\left(  i\lambda_{k}t\right)
}{C^{\prime}\left(  \lambda_{k}\right)  }\ ,\ \ \ C\left(  z\right)
=\det\left(  zI-M\right)  \ ,
\end{equation}
where I have assumed the eigenvalues are non-degenerate. \ For degenerate
eigenvalues, appropriate limits must be taken. \ 

Alternatively, the symmetric polynomials can be written as invariant traces,
hence computed directly from $M$ without knowing the individual eigenvalues
\cite{CF}.%
\begin{equation}
S_{m}\left[  \lambda\right]  =\frac{1}{m!}~\mathcal{I}_{m}\ ,
\end{equation}%
\begin{equation}
\mathcal{I}_{m}=\det\left(
\begin{array}
[c]{ccccccc}%
\operatorname*{tr}\left(  M\right)  & m-1 & 0 & \cdots & 0 & 0 & 0\\
\operatorname*{tr}\left(  M^{2}\right)  & \operatorname*{tr}\left(  M\right)
& m-2 & \cdots & 0 & 0 & 0\\
\operatorname*{tr}\left(  M^{3}\right)  & \operatorname*{tr}\left(
M^{2}\right)  & \operatorname*{tr}\left(  M\right)  & \cdots & 0 & 0 & 0\\
\vdots & \vdots & \vdots & \ddots & \vdots & \vdots & \vdots\\
\operatorname*{tr}\left(  M^{m-2}\right)  & \operatorname*{tr}\left(
M^{m-3}\right)  & \operatorname*{tr}\left(  M^{m-4}\right)  & \cdots &
\operatorname*{tr}\left(  M\right)  & 2 & 0\\
\operatorname*{tr}\left(  M^{m-1}\right)  & \operatorname*{tr}\left(
M^{m-2}\right)  & \operatorname*{tr}\left(  M^{m-3}\right)  & \cdots &
\operatorname*{tr}\left(  M^{2}\right)  & \operatorname*{tr}\left(  M\right)
& 1\\
\operatorname*{tr}\left(  M^{m}\right)  & \operatorname*{tr}\left(
M^{m-1}\right)  & \operatorname*{tr}\left(  M^{m-2}\right)  & \cdots &
\operatorname*{tr}\left(  M^{3}\right)  & \operatorname*{tr}\left(
M^{2}\right)  & \operatorname*{tr}\left(  M\right)
\end{array}
\right)  \ . \label{Invariants}%
\end{equation}
This form of the invariants follows from the generating function for the
polynomials as given by a determinant,%
\begin{equation}
\det\left(  I+tM\right)  =\sum_{m=0}^{N}t^{m}S_{m}\ .
\end{equation}

The result (\ref{ExponentialPoly}) may also be derived in a straightforward
way just by taking \href{https://en.wikipedia.org/wiki/Laplace_transform}{the
inverse Laplace transform} of the resolvent for a general complex $N\times N$
matrix $M$, as expressed by
\begin{equation}
\frac{1}{I-sM}=\sum_{n=0}^{N-1}M^{n}R_{n}\left(  s\right)  \ ,\ \ \ R_{n}%
\left(  s\right)  =\frac{s^{n}}{\det\left(  1-sM\right)  }%
~\operatorname*{Trunc}_{N-1-n}\left[  \det\left(  1-sM\right)  \right]  \ ,
\end{equation}
where $\operatorname*{Trunc}\limits_{k}\left[  f\left(  s\right)  \right]
\equiv\sum_{m=0}^{k}\frac{1}{m!}f^{\left(  m\right)  }\left(  0\right)  s^{m}$
is the $k$th-order Taylor polynomial for the function $f$. \ This result for
the resolvent is well-known \cite{He,TSvK,TLC}. \ In this approach the
response function is encountered as
\begin{equation}
F\left(  t\right)  =%
{\displaystyle\oint}
\frac{e^{itz}}{C\left(  z\right)  }~dz\ ,
\end{equation}
where the counter-clockwise contour integral $%
{\displaystyle\oint}
$ encloses all the eigenvalues of $M$.

To be more specific I now give some results appropriate to the group
$SU\left(  N\right)  $, for various $N$, that evince most of the general
features. \ In these specific cases $I$ is always the $N\times N$ unit matrix
and $H$ is a \emph{traceless, hermitian matrix}, so the eigenvalues are all
\emph{real}. \ In addition, $F_{N}\left(  t\right)  $ is the response
function\ for the $N\times N$ matrix $H$, as defined above.\ The corresponding
characteristic function and its derivative are%
\begin{equation}
C\left(  z\right)  =%
{\displaystyle\prod\limits_{k=1}^{N}}
\left(  z-\lambda_{k}\right)  \ ,\ \ \ C^{\prime}\left(  z\right)  =\sum
_{k=1}^{N}\left(
{\displaystyle\prod\limits_{\substack{m=1\\m\neq k}}^{N}}
\left(  z-\lambda_{m}\right)  \right)  \ ,\ \ \ C^{\prime}\left(  \lambda
_{k}\right)  =%
{\displaystyle\prod\limits_{\substack{m=1\\m\neq k}}^{N}}
\left(  \lambda_{k}-\lambda_{m}\right)  \ .
\end{equation}
Again, I will assume the eigenvalues are non-degenerate. \ Otherwise, for the
expressions to follow appropriate limits must be taken. \ 

For the $SU\left(  2\right)  $ doublet, as is well-known,%
\begin{equation}
\exp\left(  itH\right)  =\left[  H-iI\frac{d}{dt}\right]  F_{2}\left(
t\right)  \ .
\end{equation}
For the $SU\left(  3\right)  $ triplet, as given in \cite{CZ},%
\begin{equation}
\exp\left(  itH\right)  =\left[  H^{2}-iH\frac{d}{dt}-I\left(  \tfrac{1}%
{2}\operatorname*{tr}\left(  H^{2}\right)  +\frac{d^{2}}{dt^{2}}\right)
\right]  F_{3}\left(  t\right)  \ .
\end{equation}
For the $SU\left(  4\right)  $ quartet,%
\begin{equation}
\exp\left(  itH\right)  =\left[  H^{3}-iH^{2}\frac{d}{dt}-H\left(  \tfrac
{1}{2}\operatorname*{tr}\left(  H^{2}\right)  +\frac{d^{2}}{dt^{2}}\right)
+I\left(  -\tfrac{1}{3}\operatorname*{tr}\left(  H^{3}\right)  +\tfrac{1}%
{2}i\operatorname*{tr}\left(  H^{2}\right)  \frac{d}{dt}+i\frac{d^{3}}{dt^{3}%
}\right)  \right]  F_{4}\left(  t\right)  \ .
\end{equation}
And finally, for the $SU\left(  5\right)  $ quintet,%
\begin{align}
\exp\left(  itH\right)   &  =\left[  H^{4}-iH^{3}\frac{d}{dt}-H^{2}\left(
\tfrac{1}{2}\operatorname*{tr}\left(  H^{2}\right)  +\frac{d^{2}}{dt^{2}%
}\right)  +H\left(  -\tfrac{1}{3}\operatorname*{tr}\left(  H^{3}\right)
+\tfrac{1}{2}i\operatorname*{tr}\left(  H^{2}\right)  \frac{d}{dt}%
+i\frac{d^{3}}{dt^{3}}\right)  \right. \nonumber\\
&  \left.  +I\left(  \tfrac{1}{8}\left(  \operatorname*{tr}\left(
H^{2}\right)  \right)  ^{2}-\tfrac{1}{4}\operatorname*{tr}\left(
H^{4}\right)  +\tfrac{1}{3}\operatorname*{tr}\left(  H^{3}\right)  i\frac
{d}{dt}+\tfrac{1}{2}\operatorname*{tr}\left(  H^{2}\right)  \frac{d^{2}%
}{dt^{2}}+\frac{d^{4}}{dt^{4}}\right)  \right]  F_{5}\left(  t\right)  \ .
\end{align}
Etc. \ Note how the form of the result for $SU\left(  N\right)  $ can be
immediately obtained from that for $SU\left(  N+1\right)  $ by discarding the
unit matrix term for the latter, by decrementing the exponents of the
remaining matrix powers of $H$ (but do \emph{not} change the exponents inside
invariant traces), and by replacing $F_{N+1}\rightarrow F_{N}$. \ 

Conversely, the result for $SU\left(  N\right)  $ can be obtained from that
for $SU\left(  N-1\right)  $ just by replacing $F_{N-1}\rightarrow F_{N}$, by
incrementing the powers of $H$ in the $SU\left(  N-1\right)  $ expression, and
finally by adding the appropriate unit matrix term. \ That unit matrix term
for $SU\left(  N\right)  $ can always be expressed as a series of derivatives
of $F_{N}$ with \textquotedblleft Vi\`{e}te coefficients\textquotedblright%
\ given by invariant traces, as described above, namely,%
\begin{equation}
I\times\left(  -1\right)  ^{N-1}\sum_{n=0}^{N-1}\frac{1}{n!}~\mathcal{I}%
_{n}\left(  i\frac{d}{dt}\right)  ^{N-1-n}F_{N}\left(  t\right)  \ .
\label{UnitMatrixTerm}%
\end{equation}
The sum here can also be written as an inverse Laplace transform:%
\begin{equation}
\mathcal{L}^{-1}\left[  1-\left(  -is\right)  ^{N}\frac{\det\left(  H\right)
}{\det\left(  I-isH\right)  }\right]  =\left(  -1\right)  ^{N-1}\sum
_{n=0}^{N-1}\frac{1}{n!}~\mathcal{I}_{n}\left(  i\frac{d}{dt}\right)
^{N-1-n}F_{N}\left(  t\right)  \ . \label{InverseLaplace}%
\end{equation}
The explicit form of the $\mathcal{I}_{k}$ invariants for all $k$ are given in
(\ref{Invariants}) as traced powers of $H$. \ For the specific examples it
suffices to note: \ $\mathcal{I}_{0}=1$, $\mathcal{I}_{1}=\operatorname*{tr}%
\left(  H\right)  $, $\mathcal{I}_{2}=\left(  \operatorname*{tr}\left(
H\right)  \right)  ^{2}-\operatorname*{tr}\left(  H^{2}\right)  $,
$\mathcal{I}_{3}=\left(  \operatorname*{tr}\left(  H\right)  \right)
^{3}-3\operatorname*{tr}\left(  H\right)  \operatorname*{tr}\left(
H^{2}\right)  +2\operatorname*{tr}\left(  H^{3}\right)  $, and $\mathcal{I}%
_{4}=\left(  \operatorname*{tr}\left(  H\right)  \right)  ^{4}-6\left(
\operatorname*{tr}\left(  H\right)  \right)  ^{2}\operatorname*{tr}\left(
H^{2}\right)  +8\operatorname*{tr}\left(  H\right)  \operatorname*{tr}\left(
H^{3}\right)  +3\left(  \operatorname*{tr}\left(  H^{2}\right)  \right)
^{2}-6\operatorname*{tr}\left(  H^{4}\right)  $, where $\operatorname*{tr}%
\left(  H\right)  =0$ for $SU\left(  N\right)  $ generators.

Thus a compact \emph{form} of the fundamental exponential polynomial for any
$SU\left(  N\right)  $ can be obtained by using (\ref{UnitMatrixTerm}) to
construct sequentially the hierarchy of polynomials for $SU\left(  M\leq
N\right)  $, starting from the trivial result $I$ for $SU\left(  1\right)  $
\cite{ThisGeneralizes}. \ However, more work is required to obtain explicit
results for the response functions since the eigenvalues for a generic
$N\times N$ matrix, whether hermitian or not, pose some challenges \cite{Pan},
especially in the limit of large $N$.

Nevertheless, as an explicit example, an element of an $SU\left(  2\right)  $
subgroup of $SU\left(  N\right)  $ is always manageable in closed form, where
$N=2j+1$ fixes the spin of the embedded $SU\left(  2\right)  $ representation.
\ The standard choice for the generator $\hat{n}\cdot\overrightarrow{J}$ of a
rotation about axis $\hat{n}$ gives a characteristic function
\begin{equation}
C\left(  \lambda\right)  =\det\left(  \lambda I-\hat{n}\cdot\overrightarrow{J}%
\right)  =%
{\displaystyle\prod\limits_{k=0}^{2j}}
\left(  \lambda-\left(  j-k\right)  \right)  =\frac{\Gamma\left(
\lambda+j+1\right)  }{\Gamma\left(  \lambda-j\right)  }\ .
\end{equation}
The trace norm of this standard generator is given by a product of the
quadratic $su\left(  2\right)  $ Casimir, $j\left(  j+1\right)  $, and the
matrix rank, $N=2j+1$, namely,
\begin{equation}
\operatorname*{tr}\left[  \left(  \hat{n}\cdot\overrightarrow{J}\right)
^{2}\right]  =%
{\displaystyle\sum\limits_{k=0}^{2j}}
\left(  j-k\right)  ^{2}=\frac{1}{3}~j\left(  j+1\right)  \times\left(
2j+1\right)  \ .
\end{equation}
Traces of all odd powers of $\hat{n}\cdot\overrightarrow{J}\ $vanish, of
course, since the eigenvalues are symmetrically distributed about zero.
$\ $Traces of all even powers are easily computed, as given by the following
generating function (i.e. the character of the group element, upon replacing
$x\rightarrow i\theta$ where $\theta$ is the rotation angle):%
\begin{gather}
\sum_{n=0}^{\infty}\frac{x^{n}}{n!}\operatorname*{tr}\left[  \left(  \hat
{n}\cdot\overrightarrow{J}\right)  ^{n}\right]  =%
{\displaystyle\sum\limits_{k=0}^{2j}}
e^{x\left(  j-k\right)  }=\frac{\sinh\left(  \left(  2j+1\right)  x/2\right)
}{\sinh\left(  x/2\right)  }\\
=\left(  2j+1\right)  \left(  1+\frac{x^{2}}{6}~j\left(  j+1\right)
+\frac{x^{4}}{360}~j\left(  j+1\right)  \times\left(  3j\left(  j+1\right)
-1\right)  +O\left(  x^{6}\right)  \right)  \ .\nonumber
\end{gather}
Note that $\frac{1}{2j+1}\operatorname*{tr}\left[  \left(  \hat{n}%
\cdot\overrightarrow{J}\right)  ^{2k}\right]  $ is always a $k$th order
polynomial in the quadratic Casimir.

More importantly for the exponentiated matrix, the response function for the
embedded $SU\left(  2\right)  $ spin $j$ representation is%
\begin{equation}
F_{N=2j+1}\left(  \theta\right)  =\sum_{k=0}^{2j}\frac{\exp\left(  i\left(
j-k\right)  \theta\right)  }{C^{\prime}\left(  j-k\right)  }=\frac{\left(
2i\right)  ^{2j}}{\left(  2j\right)  !}~\sin^{2j}\left(  \theta/2\right)  \ .
\end{equation}
Derivatives of this response function are readily evaluated to obtain a matrix
polynomial expansion for $\exp\left(  i\theta\hat{n}\cdot\overrightarrow{J}%
\right)  $, although it takes quite a bit of additional work to reduce the
resulting polynomial to the compact form in \cite{CFZ,CvK}.

For generic $SU\left(  N\right)  $ generators the eigenvalues and the response
function cannot be so easily determined as was the case for the embedded
$SU\left(  2\right)  $ example. \ However, there is an elegant geometrical way
to parameterize the eigenvalues that generalizes the well-known picture
developed by Vi\`{e}te in the 16th century for cubic equations with three real
roots \cite{Viete}. \ Indeed, the discussion to follow is easily adapted to
provide a \emph{visualization }of the solution to any $N$th order equation
with only real roots, although it does not immediately give an explicit
solution for the roots in terms of traces and/or determinants.

For $SU\left(  N\right)  $ the $N$ real eigenvalues of the traceless hermitian
generator $H$ may be envisioned as the components of a Euclidean vector,
$\overrightarrow{\lambda}\in\mathbb{E}_{N}$. \ This vector of eigenvalues lies
in an $\left(  N-1\right)  $-dimensional hyperplane\ that passes through the
origin of the Euclidean space, with normal $\overrightarrow{n}\equiv\left(
1,1,\cdots,1\right)  \in\mathbb{E}_{N}$, as defined by $0=\overrightarrow{n}%
\cdot\overrightarrow{\lambda}=\sum_{k=1}^{N}\lambda_{k}=\operatorname*{tr}%
\left(  H\right)  $. \ The vector $\overrightarrow{\lambda}$ is also a point
on a sphere $S_{N-1}\subset\mathbb{E}_{N}$, with the non-vanishing radius of
the sphere defined by the matrix invariant $r^{2}\equiv\overrightarrow{\lambda
}\cdot\overrightarrow{\lambda}=\sum_{k=1}^{N}\lambda_{k}^{2}%
=\operatorname*{tr}\left(  H^{2}\right)  $. \ The intersection of this
hyperplane and $\left(  N-1\right)  $-sphere define a particular $\left(
N-2\right)  $-sphere in the eigenvalue space, $S_{N-2}\subset\mathbb{E}_{N}$.
\ The eigenvalues of $H$ therefore comprise the components of a vector
$\overrightarrow{\lambda}$ that is a point on this $S_{N-2}$, and for a given
radius $r$ the vector $\overrightarrow{\lambda}$ may be completely specified
by the standard $N-2$ spherical polar angles parameterizing $S_{N-2}$.

For example, for $SU\left(  3\right)  $ the vector of eigenvalues is a point
on the circle defined by the intersection of the aforementioned plane and
2-sphere. \ The components of $\overrightarrow{\lambda}$ are specified by the
radius of the circle and by a single angle $\theta$, as noted long ago by
Vi\`{e}te for cubic equations with three real roots. \ Thus%
\begin{equation}
\lambda_{k}=\sqrt{\tfrac{2}{3}}~r\cos\left(  \theta+2\pi k/3\right)  \text{
\ for \ }k=1,2,3\text{, \ where }\prod\limits_{k=1}^{3}\lambda_{k}=\tfrac
{1}{3\sqrt{6}}~r^{3}\cos\left(  3\theta\right)  =\det\left(  H\right)  \ .
\label{SU3Roots}%
\end{equation}
Perhaps the only surprising features apparent in (\ref{SU3Roots}) are the
facts that the radius of the circle has been rescaled \cite{Footnote1} and
that the eigenvalues may also be viewed as the projections of three points
\emph{equally spaced} on the circle. \ That is to say, the eigenvalues are
given by projecting along a particular axis the vertices of an equilateral
triangle (i.e. a 2-simplex) circumscribed by a circle, where the axis of
interest is a diameter of the circle. \ The angle $\theta$ is determined by
the orientation of the triangle with respect to this axis.\ \ 

Perhaps even more surprisingly, this geometrical picture generalizes
completely to any $N$. \ The eigenvalues of any generator for the fundamental
representation of $SU\left(  N\right)  $ are given by projecting the $N$
vertices of an $\left(  N-1\right)  $-simplex onto a particular axis that
passes through the center of the simplex. \ The orientation of the simplex
with respect to the axis determines the $N-2$ angles that are needed, in
addition to the radius $r$, to specify the components of
$\overrightarrow{\lambda}$.\ \ I now explain in detail how these geometrical
facts are established, if they are not already obvious.

The components of $\overrightarrow{\lambda}$ are the eigenvalues of $H$.
\ Each of these components is given by a Euclidean inner product with one of
the unit vectors that define the standard orthonormal basis for $\mathbb{E}%
_{N}$. \ Thus
\begin{equation}
\lambda_{k}=\widehat{e}_{k}\cdot\overrightarrow{\lambda}\ ,\ \ \ \text{where}%
\ \ \ \widehat{e}_{1}=\left(  1,0,\cdots,0\right)  \ ,\ \ \ \widehat{e}%
_{2}=\left(  0,1,\cdots,0\right)  \ ,\ \cdots\ \ \ ,\ \ \ \widehat{e}%
_{N}=\left(  0,\cdots,0,1\right)  \ .
\end{equation}
On the other hand, these $N$ unit vectors are the vertices of the $N$-cell
that defines
\href{https://en.wikipedia.org/wiki/Simplex#The_standard_simplex}{the standard
$\left(  N-1\right)  $-simplex}, embedded in $\mathbb{E}_{N}$. \ This standard
simplex lies in a hyperplane \textquotedblleft parallel\textquotedblright\ to
the one that contains the vector of eigenvalues, with a normal again given by
the aforementioned $\overrightarrow{n}$. \ Moreover, since $\overrightarrow{n}%
\cdot\overrightarrow{\lambda}=0$, the inner products with
$\overrightarrow{\lambda}$\ are unchanged if each of the $\widehat{e}_{k}$ is
translated by adding any amount of the normal $\overrightarrow{n}$. \ Thus we
may rigidly translate the standard simplex so that it too is centered on the
origin and lies in the hyperplane containing the eigenvalues. \ Upon doing so,
we obtain a simplex with vertices $\widehat{e}_{k}-\frac{1}{N}%
\overrightarrow{n}$, where the eigenvalues are now given by inner products
$\lambda_{k}=\left(  \widehat{e}_{k}-\frac{1}{N}\overrightarrow{n}\right)
\cdot\overrightarrow{\lambda}$. \ More explicitly,%
\begin{equation}
\widehat{e}_{1}-\tfrac{1}{N}\overrightarrow{n}=\left(  \tfrac{N-1}{N}%
,-\tfrac{1}{N},\cdots,-\tfrac{1}{N}\right)  \ ,\ \ \ \widehat{e}_{2}-\tfrac
{1}{N}\overrightarrow{n}=\left(  -\tfrac{1}{N},\tfrac{N-1}{N},\cdots
,-\tfrac{1}{N}\right)  \ ,\ \cdots\ \ \ ,\ \ \ \widehat{e}_{N}-\tfrac{1}%
{N}\overrightarrow{n}=\left(  -\tfrac{1}{N},\cdots,-\tfrac{1}{N},\tfrac
{N-1}{N}\right)  \ .
\end{equation}
In addition, we may rescale these vertices to obtain a simplex that also lies
on the $\left(  N-2\right)  $-sphere that contains the vector of eigenvalues.
\ Thus take
\begin{equation}
\overrightarrow{f}_{k}=\left(  \widehat{e}_{k}-\frac{1}{N}\overrightarrow{n}%
\right)  r~\sqrt{\frac{N}{N-1}}\ .
\end{equation}
One readily checks that these $N$ vectors are indeed the vertices of a
simplex, with
\begin{equation}
\overrightarrow{f}_{k}\cdot\overrightarrow{f}_{k}=r^{2}\text{ \ for any
\ }k=1,2,\cdots,N\text{ \ and\ \ }\overrightarrow{f}_{k}\cdot
\overrightarrow{f}_{m}=\frac{1}{1-N}~r^{2}\text{ \ for any \ }k\neq m\ .
\end{equation}
This simplex corresponds to the array of weight vectors \cite{Wybourne} for
the fundamental representation of $SU\left(  N\right)  $. \ The eigenvalues
are now obtained by projecting the vertices of this final $\left(  N-1\right)
$-simplex onto an axis defined by the direction of $\overrightarrow{\lambda}$,
i.e. by projecting the vertices onto a particular diameter of\ the $\left(
N-2\right)  $-sphere.

Thus, with unit vector $\widehat{e}=\overrightarrow{\lambda}/r$,%
\begin{equation}
\lambda_{k}=\sqrt{\frac{N-1}{N}}~\overrightarrow{f}_{k}\cdot\widehat{e}\ .
\label{SUnRoots}%
\end{equation}
But note, what matters for the eigenvalues is the relative orientation of the
simplex with respect to the axis. \ Thus we may produce the same eigenvalues
if we choose the axis to be in any convenient direction, and appropriately
orient an equivalent simplex with respect to the chosen axis.

For example, for the eigenvalues of the $SU\left(  3\right)  $ generator, as
given in (\ref{SU3Roots}), we may choose the axis as the horizontal on the
plane, with $\theta$ the planar angle measured counterclockwise from the
horizontal. \ We may then visualize the eigenvalues as projections onto the
horizontal axis of the vertices of an equilateral triangle in the plane --- a
2-simplex --- those vertices being located at angles $\theta$ and $\theta
\pm2\pi/3$, and at a distance $\sqrt{\frac{2}{3}}~r$ from the origin.

For other examples, consider $SU\left(  4\right)  $ and $SU\left(  5\right)
$, and note that the freedom to re-orient jointly the axis \emph{and} the
simplex, as if they comprised a single rigid body, leads to apparently
different --- but nevertheless equivalent --- parameterizations of the
eigenvalues in terms of the standard spherical polar angles for $S_{N-2}$. \ 

For $SU\left(  4\right)  $ the vector of eigenvalues is a point on a 2-sphere.
\ The four eigenvalues, i.e. the components of the four-vector
$\overrightarrow{\lambda}$,\ may be specified in terms of two polar angles.%
\begin{align}
\lambda_{1}  &  =r\left(  -\frac{1}{\sqrt{2}}\sin\phi\sin\theta-\frac{1}%
{2}\cos\theta\right)  \ ,\\
\lambda_{2}  &  =r\left(  +\frac{1}{\sqrt{2}}\sin\phi\sin\theta-\frac{1}%
{2}\cos\theta\right)  \ ,\\
\lambda_{3}  &  =r\left(  -\frac{1}{\sqrt{2}}\cos\phi\sin\theta+\frac{1}%
{2}\cos\theta\right)  \ ,\\
\lambda_{4}  &  =r\left(  +\frac{1}{\sqrt{2}}\cos\phi\sin\theta+\frac{1}%
{2}\cos\theta\right)  \ .
\end{align}
These are just projections of the vertices of
\href{https://en.wikipedia.org/wiki/Tetrahedron}{a regular tetrahedron} onto
an axis through its center, whose orientation is left as an exercise for the
reader \cite{Footnote2}, but which is easily visualized in three Euclidean
dimensions.%
\begin{center}
\includegraphics[
height=3.2091in,
width=4.8003in
]%
{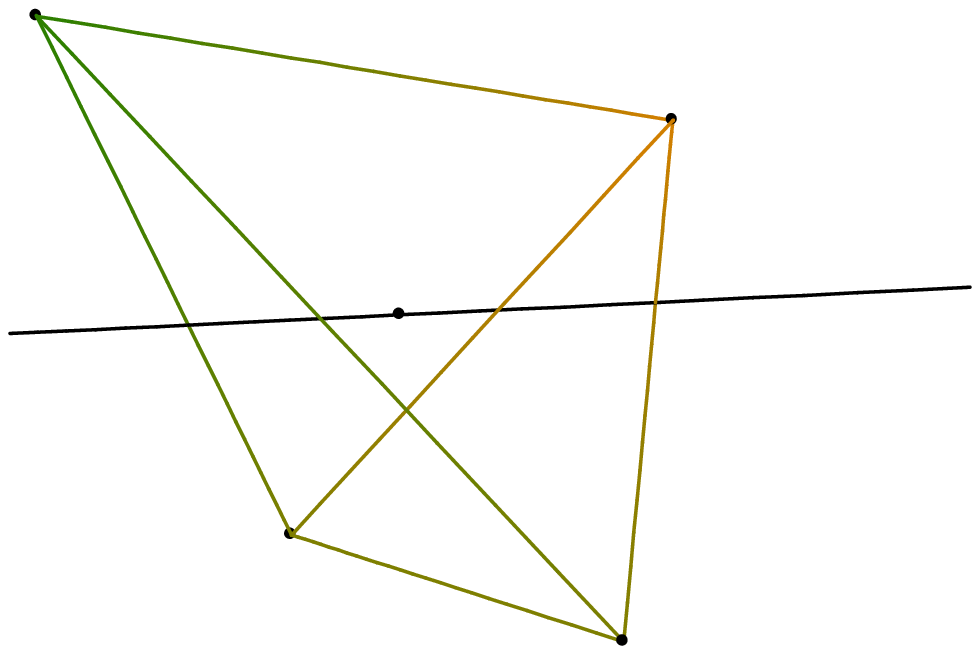}%
\\
A tetrahedron pierced by an axis through its center.
\end{center}

For $SU\left(  5\right)  $ the vector of eigenvalues is a point on a 3-sphere.
\ The five eigenvalues, i.e. the components of the five-vector
$\overrightarrow{\lambda}$,\ may be specified in terms of three polar angles.%

\begin{align}
\lambda_{1}  &  =r\left(  -\frac{4}{2\sqrt{5}}\cos\psi\right)  \ ,\\
\lambda_{2}  &  =r\left(  \frac{1}{2\sqrt{5}}\cos\psi-\frac{3}{2\sqrt{3}}%
\cos\theta\sin\psi\right)  \ ,\\
\lambda_{3}  &  =r\left(  \frac{1}{2\sqrt{5}}\cos\psi+\frac{1}{2\sqrt{3}}%
\cos\theta\sin\psi-\frac{2}{\sqrt{6}}\cos\phi\sin\theta\sin\psi\right)  \ ,\\
\lambda_{4}  &  =r\left(  \frac{1}{2\sqrt{5}}\cos\psi+\frac{1}{2\sqrt{3}}%
\cos\theta\sin\psi+\frac{1}{\sqrt{6}}\cos\phi\sin\theta\sin\psi-\frac{1}%
{\sqrt{2}}\sin\theta\sin\phi\sin\psi\right)  \ ,\\
\lambda_{5}  &  =r\left(  \frac{1}{2\sqrt{5}}\cos\psi+\frac{1}{2\sqrt{3}}%
\cos\theta\sin\psi+\frac{1}{\sqrt{6}}\cos\phi\sin\theta\sin\psi+\frac{1}%
{\sqrt{2}}\sin\theta\sin\phi\sin\psi\right)  \ .
\end{align}
These are just projections of the vertices of a
\href{https://en.wikipedia.org/wiki/5-cell}{pentatope} onto an axis through
its center. \ In this case I have chosen things so that the relative
orientations of the axis and the five vertices are easily determined in terms
of the angular parameterization, even if the entire geometry is not so easily
visualized as in the $SU\left(  4\right)  $ example.

Finally, just as the eigenvalues are readily expressed in terms of the radius
and the various angles, so too are the latter expressed in terms of
combinations of the eigenvalues. \ But it is not easy --- there is no royal
road to find the roots of quintic and higher order equations --- to express
the eigenvalues, or the corresponding angles, in terms of the invariant traces
or the determinant of $H$ for $SU\left(  N>4\right)  $. \ 

For example, for $SU\left(  4\right)  $ as parameterized above,
\begin{align}
\operatorname*{tr}\left(  H\right)   &  =\sum_{k=1}^{4}\lambda_{k}%
=0\ ,\ \ \ \operatorname*{tr}\left(  H^{2}\right)  =\sum_{k=1}^{4}\lambda
_{k}^{2}=r^{2}\ ,\\
\operatorname*{tr}\left(  H^{3}\right)   &  =\sum_{k=1}^{4}\lambda_{k}%
^{3}=\tfrac{3}{4}r^{3}\sin\theta\sin\left(  2\theta\right)  \cos\left(
2\phi\right)  \ ,\\
\det\left(  H\right)   &  =%
{\displaystyle\prod\limits_{k=1}^{4}}
\lambda_{k}=\tfrac{1}{8}\left(  \left(  \operatorname*{tr}\left(
H^{2}\right)  \right)  ^{2}-2\operatorname*{tr}\left(  H^{4}\right)  \right)
\nonumber\\
&  =\tfrac{1}{16}r^{4}\left(  1+\left(  2\sin^{2}\phi-3\right)  \sin^{2}%
\theta\right)  \left(  1+\left(  2\cos^{2}\phi-3\right)  \sin^{2}%
\theta\right)  \ .
\end{align}
This leads to equations cubic and quadratic in $\sin^{2}\theta$ and $\sin
^{2}\phi$ that can be solved in closed form in terms of elementary functions
of the traces of $H^{2}$, $H^{3}$, and $H^{4}$ or $\det\left(  H\right)  $.
\ For detailed discussions on how to solve the resulting equations, see
\cite{Kusnezov,Quartic}.

For $SU\left(  5\right)  $ as parameterized above, it is straightforward to
check that%
\begin{equation}
\operatorname*{tr}\left(  H^{3}\right)  =r^{3}\left(  \frac{3}{\sqrt{5}%
}\left(  \cos\psi\right)  \left(  \frac{1}{2}-\cos^{2}\psi\right)  +\frac
{5}{2\sqrt{3}}\left(  \sin^{3}\psi\cos\theta\right)  \left(  \frac{3}{5}%
-\cos^{2}\theta\right)  +\frac{2\sqrt{2}}{\sqrt{3}}\left(  \sin^{3}\psi
\sin^{3}\theta\cos\phi\right)  \left(  \frac{3}{4}-\cos^{2}\phi\right)
\right)  \ ,
\end{equation}
etc. \ However, in this case I can not find an \emph{elementary} expression of
the angles, or eigenvalues, in terms of invariant traces of powers of $H$.
\ But then, neither can anyone else, as is well-known \cite{Quintic&Beyond}.

Nonetheless, at the very least I hope the constructions described here for any
$N$ may serve as elementary examples that illustrate the elegant interplay
between group theory and geometry. \ As H S M Coxeter once said \cite{Coxeter},

\begin{quote}
{\small \textquotedblleft\ ... in four or more dimensions, we can never fully
comprehend them by direct observation. \ In attempting to do so, however, we
seem to peep through a chink in the wall of our physical limitations, into a
new world of dazzling beauty. \ Such an escape from the turbulence of ordinary
life will perhaps help to keep us sane.\textquotedblright}
\end{quote}

\paragraph{Acknowledgement}

I thank Professor Curtright for focusing my attention on \cite{CZ}, and for
suggesting that the simple geometrical construction so succinctly described
therein would generalize to any $SU\left(  N\right)  $.


\begin{thebibliography}{99}
\bibitem {ThisGeneralizes}The same procedure may be used to construct
sequentially a hierarchy of polynomials for the exponentials\ of almost any
matrix, as in (\ref{ExponentialPoly}).
\bibitem {Footnote1}Although the eigenvalues themselves still define a circle
with the same radius as before, namely $r^{2}$, since $\sum_{k=1}^{3}\cos
^{2}\left(  \theta+2\pi k/3\right)  =3/2$.
\bibitem {Footnote2}Here is the solution to the exercise: \ The tetrahedron's
vertices are at a distance $\sqrt{3}~r/2$ from the origin, such that
$\lambda_{1}=\frac{1}{2}r\left(  -1,1,-1\right)  \cdot\widehat{e}$,
$\lambda_{2}=\frac{1}{2}r\left(  -1,-1,1\right)  \cdot\widehat{e}$,
$\lambda_{3}=\frac{1}{2}r\left(  1,-1,-1\right)  \cdot\widehat{e}$, and
$\lambda_{4}=\frac{1}{2}r\left(  1,1,1\right)  \cdot\widehat{e}$, where
$\widehat{e}\equiv\left(  \cos\theta,\sin\theta\cos\left(  \phi+\pi/4\right)
,\sin\theta\sin\left(  \phi+\pi/4\right)  \right)  $.
\end{thebibliography}
\end{document}